\documentclass[a4paper,12pt]{article}
\usepackage{amsmath,amssymb,amsfonts,amsthm}
\usepackage{amscd}
\newtheorem{theorem}{Theorem}[section]
\newtheorem{prop}[theorem]{Proposition}
\newtheorem{lemma}[theorem]{Lemma}
\newtheorem{coro}[theorem]{Corollary}
\newtheorem{exa}[theorem]{Example}
\newtheorem{defi}[theorem]{Definition}
\newtheorem{rem}[theorem]{Remark}
\newcommand{\bte}{\begin{theorem}\quad  }
\newcommand{\ete}{\end{theorem} }
\newcommand{\bpr}{\begin{prop}\quad  }
\newcommand{\epr}{\end{prop} }
\newcommand{\ble}{\begin{lemma}\quad }
\newcommand{\ele}{\end{lemma}}
\newcommand{\bco}{\begin{coro}\quad }
\newcommand{\eco}{\end{coro} }
\newcommand{\bex}{\begin{exa}\quad \rm }
\newcommand{\eex}{\end{exa} }
\newcommand{\bde}{\begin{defi}\quad \rm }
\newcommand{\ede}{\end{defi} }
\newcommand{\brm}{\begin{rem} \quad \rm}
\newcommand{\erm}{\end{rem} }
\newcommand{\bpf}{\begin{proof}[\bf{Proof.\quad}] \rm}
\newcommand{\epf}{ \end{proof}}
\newcommand{\bdm}{\begin{displaymath} }
\newcommand{\edm}{\end{displaymath} }

\newcommand{\lb}{\label}
\newcommand{\lo}{\longrightarrow}

\begin{document}
\title{{\textbf{ Co-uniform and hollow $S$-acts over monoids }}}
\author{  \textbf{R.  Khosravi\footnote{Department of Mathematics, Faculty of Sciences, Fasa University, Fasa, Iran.}  .   M. Roueentan\footnote{ College of Engineering, Lamerd Higher Education Center, Lamerd, Iran  } 
}}

\date{
}
\maketitle
\begin{abstract}

In this paper, we first introduce the notions of superfluous and coessential subacts. Then    hollow and co-uniform $S$-acts are defined as the acts that all proper subacts are   superfluous  and coessential, respectively.  Also it is indicated that the class of hollow $S$-acts is properly between two classes of indecomposable and locally cyclic $S$-acts. Moreover, using the notion of radical of an $S$-act as  the intersection of all maximal subact, the relations between   hollow and local $S$-acts are investigated. Ultimately, the notion of a  supplement of a subact is defined to characterize the union of   hollow $S$-act.
\vspace{0.25 cm} \textbf{\\Key Words} ~  monoids . $S$-acts .
 superfluous . coessential .   hollow
  \\\textbf{AMS 2010 Mathematics subject
classification} 06F05 . 20M30
\end{abstract}\footnote{ $\boxtimes$ Roghaieh Khosravi (Corresponding Author)

 khosravi@fasau.ac.ir

  Mohammad Roueentan

 rooeintan@lamerdhec.ac.ir
 
 
}
{\tiny{\tiny{\tiny{\section{\large{\hspace{-0.7 cm}.\quad
Introduction}}}}}}

 A submodule $K$ of an $R$-module  $M$ is called   superfluous (small), if the equality $N + K = M$ implies that $K = M$. The notion of small submodule plays a fundamental role in the category of modules over rings.  According to \cite{flu74}, a non-zero module $M$ is defined hollow if every submodule  of $M$ is  small (superfluous). The  classical notion of  hollow modules has been studied extensively for a long time in many papers (see for example \cite{hara, wis}). In the category of $S$-acts the notions of small (coessential) and superfluous subacts are distinct which we define both as follows. For $S$-acts, first we refer the reader to \cite{kilp2000} for preliminaries and basic results related monoids and $S$-acts. A subact $B_S$ of $A_S$ is called large in $A_S$ if any homomorphism $g : A_S \lo C_S$ such that $g|_B$ is a monomorphism is itself a monomorphism. An extension $B$ of $A$ with the embedding $f : A_S\lo B_S$ is called an essential extension of $A$ if $\text{Im}f$ is large in $B$.

The categorical dual of essential extension is called  coessential epimorphism which we recall  as follows.
Let $S$ be a monoid. An act $B_S$ is called a \textit{cover} of an
act $A_S$ if there exists an epimorphism $f: B_S\rightarrow A_S$
such that for any proper subact $C_S$ of $B_S$ the restriction
$f\mid_{C_S}$ is not an epimorphism. An epimorphism with this
property is called a \textit{coessential epimorphism}. Indeed it is defined in order to investigate  $X$-perfect monoids as monoids over which every right $S$-act has an  X-cover, where $X$ is an act
property which is preserved under coproduct. More  information about  various kinds of cover of acts one can see \cite{isb,kh, kilp96a, mah}.

As a dual of large subact, we call $B_S$ a coessential (small) subact of $A_S$ if $A_S$ is a cover of the Rees factor act $A_S/B_S$. According to the notion of superfluous submodule, a subact $B_S$ of an $S$-act $A_S$ shall be called superfluous if the union of $B_S$ with every proper subact of $A_S$ is also a proper subact of $A_S$.  In Section 2, We consider the properties of coessential and superfluous subacts. In \cite{uniform}, the authors investigated uniform acts over a semigroup $S$, as  $S$-acts that  all their  non-zero subacts are large. In module theory, the dual notion of a uniform module is that of a  hollow module. In fact  hollow and co-uniform modules are equal. For $S$-acts, as we mentioned earlier, the notion of coessential and superfluous are distinct, so we define co-uniform as a dual of uniform $S$-acts and  hollow $S$-acts with respect to the definition of  hollow in  module theory.   In Section 3, we characterize the classes of  co-uniform  and  hollow acts as the acts all proper subacts are coessential and superfluous respectively.  In Section 4, we investigate radical of an $S$-acts and local $S$-acts, and consider the relationship between local and hollow $S$-acts.
 Finally, in Section 5, a supplement of a subact and supplemented $S$-acts are introduced and  using these notions to characterize the union of  hollow $S$-acts. 
The following lemma is clearly proved which is needed in the sequel.
\ble \lb{le1.1} If $M$ is a maximal subact of a right $S$-act $A_S$, then $A/M$ is finitely generated. \ele

{\tiny{\tiny{\tiny{\section{\large{\hspace{-0.7 cm}.\quad
Coessential or Superfluous Subacts}}}}}}
In this section we introduce the notions of coessential and superfluous subacts, and  consider  general properties of them.

\bde A subact $B_S$ of an $S$-act $A_S$ is called\begin{itemize}
	
	\item[(i)]    \textit{coessential} if the epimorphism $\pi: A_S\lo A_S/{B_S}$ is a coessential epimorphism; in other words, $A_S$ is a cover of $A_S/{B_S}$. It is denoted by $B\ll A$.

\item[(ii)]    \textit{superfluous} if $B_S\cup C_S\neq A_S$ for each proper subact $C_S$ of $A_S$, and it is denoted by $B\leq_s A$.\end{itemize} \ede

In the following lemma we present an equivalent condition for being coessential.
\ble   A subact $B_S$ of an $S$-act $A_S$ is coessential if and only if for each proper subact $C_S$ of $A_S$, $C\cap B\neq \emptyset$ implies that $C\cup B\neq A$. \ele
\bpf Necessity.  Let $C_S$ be a subact of $A_S$ and $C\cap B\neq \emptyset$. As we know, $\pi|_{C_S}$ is not an epimorphism, which implies  the existence of $a\in A_S$ such that $[a]\notin \pi(C)$. Now we claim that $a\notin  C\cup B$. Otherwise, either $a\in C$ which means $[a]\in \pi(C)$ or  $a\in B$ which implies $[a]=[b]\in \pi(C)$ for some $b\in C\cap B$. Thus $C\cup B\neq A$.

Sufficiency. Let $C_S$ be a proper subact of $A_S$. If $C\cap B= \emptyset$, clearly for each $b\in B$ we have $[b]\notin \pi(C)$. Otherwise, if $C\cap B\neq \emptyset$, by assumption $C\cup B\neq A$. So we have $[a]\notin \pi(C)$ for each $a\in A\setminus (C\cup B)$. Therefore, $\pi|_{C_S}$ is not an epimorphism. \epf

In view of the previous lemma, it is obvious that being superfluous subact implies coessential. But  the converse is not valid.
For instance, let  $S$ be an arbitrary monoid and $A_S=\Theta\coprod\Theta=\{\theta_1,\theta_2\}$. Then $\{\theta_1\}$ is coessential but not superfluous. 
\ble \lb{le2.1} A coessential subact of each indecomposable right $S$-act is superfluous.\ele

\bpf  Suppose that $B$ is a coessential subact of an indecomposable right $S$-act $A_S$ and $B\cup C=A$ for a subact $C$ of $A$. If $B\cap C=\emptyset$, then $A=B\coprod C$ which contradicts with  being indecomposable. So $B\cap C\neq\emptyset$ and  $B\cup C=A$ which imply that $C=A$. Therefore, $B$ is superfluous.\epf

\ble \lb{le2.2} Suppose that $A_S,~B_S,~C_S, D_S$  are $S$-acts such that $D_S\subseteq C_S\subseteq B_S\subseteq A_S$. The following hold.
\begin{itemize}

\item[(i)] $B\leq_s A$ if and only if $C\leq_s A$ and  $B/C\leq_s A/C$.

\item[(ii)] If  $C\leq_s B$, then $C\leq_s A$.

\item[(iii)] $B\leq_s A$ if and only if for each $S$-act $X_S$ and $h: X\lo A$, $\rm{Im}h\cup B= A$ implies $\rm{Im}h=A$.
\item[(iv)] $B/D\leq_sA/D$ if and only if $B/C\leq_s A/C$ and $C/D \leq_s A/D$.

\end{itemize}
\ele
\bpf (i). Necessity. The first part is obvious. Let $K$ be a subact of $A/C$ with $ B/C\cup K=A/C$. So $D=\{t\in A|~[t]\in B/C\}$ is a subact of $A_S$ and it is easily checked that  $D\cup B=  A$. By assumption, $D= A$, and thus $K=A/C$.

Sufficiency. Let $D$ be a subact of $A$ and $D\cup B=  A$. So $ B/C\cup (D\cup C)/C=A/C$  which implies $(D\cup C)/C=A/C$. Then $D\cup C=A$ implies that $D=A$, as desired.

Parts (ii) and (iii) are clear. 

(iv)  We only show the sufficiency. Suppose that $(B/D)\cup K=A/D$ for some subact $K$ of $A/D$. Get $X=\{t\in A|~[t]\in K\}$ which is clearly a subact of $A_S$. Then $(B/C)\cup((X\cup C)/C) =A/C$. Since $B/C\leq_s A/C$, we have $X\cup C=A$. So $(C/D)\cup K =A/D$ and since $C/D \leq_s A/D$, $K=A/D$. Therefore  $B/D\leq_s A/D$.
  \epf
Similar to the proof of the previous lemma, two following lemmas are easily checked.

\ble  The following hold for a monoid $S$.
\begin{itemize}

\item[(i)] If  $C_S\subseteq B_S\subseteq A_S$ and $C\ll B$, then $C\ll A$.

\item[(ii)] If $C_S\subseteq B_S\subseteq A_S$ and $B\ll A$, then $C\ll A$ and  $B/C\ll A/C$.

\item[(ii)] If $B\ll  A$ ($B\leq_s A$) and $f:A\lo C$ is a monomorphism, then  $f(B)\ll  C$ ($f(B)\leq_s C$).

 \end{itemize} \ele
\ble  \lb{le2.3} Let $B,C$ are proper subacts of $A_S$. Then $B\cup C\leq_s A$ if and only if $B\leq_s A$ and  $C\leq_s A$.\ele

\ble\lb{le1.01}  Suppose that $B_i$ is a proper subact of $A_i$  for each $i\in I$. The following hold for a monoid $S$.
\begin{itemize}

\item[(i)]    $\coprod_{i\in I} B_i\leq_s \coprod_{i\in I} A_i$ if and only if $B_i\leq_s A_i$  for each $i\in I$.
\item[(ii)] If   $\coprod_{i\in I} B_i\ll\coprod_{i\in I} A_i$, then $B_i\ll A_i$ for each $i\in I$.

    \item[(iii)]     If  $B_i\leq_s A_i$  ($B_i\ll A_i$) for each $i\in \{1,...,n\}$, then $\cup_{i=1}^{i=n} B_i\leq_s \cup_{i=1}^{i=n} A_i$   ($\cup_{i=1}^{i=n} B_i\ll \cup_{i=1}^{i=n} A_i)$.

  \end{itemize}
     \ele

     \bpf  (i). Necessity. Suppose that   $\coprod_{i\in I} B_i\leq_s \coprod_{i\in I} A_i$. Fix $j\in I$ and $D_j$ a subact of $A_j$ such that $B_j\cup D_j=A_j$. Then $D= (\coprod_{i\neq j} A_i)\coprod D_j$ is a subact of $\coprod_{i\in I} A_i$ and $\coprod_{i\in I} B_i \cup D= \coprod_{i\in I} A_i$. By assumption, $D=\coprod_{i\in I} A_i$ which implies that $D_j=A_j$.

     Sufficiency. Suppose that $B_i\leq_s A_i$  for each $i\in I$.  Let $D$ be a subact of $\coprod_{i\in I} A_i$ such that $\coprod_{i\in I} B_i \cup D= \coprod_{i\in I} A_i$. Since $B_i$ is proper subact of $A_i$  for each $i\in I$, $D=\coprod_{i\in I} D_i$ such that $D_i\neq \emptyset$ is a subact of $A_i$. Obviously, $B_i\cup D_i=A_i$ for every $i\in I$ and by assumption $D_i=A_i$ which gives that $D= \coprod_{i\in I} A_i$.

     By a similar argument one can prove part (ii). Part (iii) is a straightforward consequence of Lemmas \ref{le2.2} and \ref{le2.3}.\epf


{\tiny{\tiny{\tiny{\section{\large{\hspace{-0.7 cm}.\quad
Co-uniform and  Hollow $S$-acts}}}}}}
In this section we study the classes of co-uniform and hollow $S$-acts.

\bde An $S$-act $A_S$ is called \textit{ co-uniform} if all proper subacts of $A_S$ are coessential, and $A_S$ is said to be \textit{ hollow} if every its proper subact  is superfluous.    \ede

Obviously,  hollow implies  co-uniform, but the converse is not valid.
 Let  $S$ be an arbitrary monoid. It is easily checked that, $\Theta\coprod\Theta$ is co-uniform but not  hollow.

\bpr Every factor act of a  (co-uniform) hollow act is also (co-uniform) hollow.  \epr

\bpf Let $A$ be a  hollow $S$-act and $f:A\lo C$  an epimorphism. Let $D$ be a proper subact of $C$. We show that $D\leq_sC$. Clearly, $B=f^{-1}(D)$ is also a proper subact of $A$. So $B\leq_sA$. Now, suppose that $D\cup E=C$. It is easily checked that $B\cup f^{-1}(E)=A$. So by assumption, $f^{-1}(E)=A$, and thus $E=C$. By a similar argument one could prove for co-uniform acts.\epf




Recall that an $S$-act $A_S$ is called locally cyclic if for all $a,a'\in A_S$ there exists $a''\in A$ such that $a,a'\in a''S$. Every locally cyclic $S$-act is indecomposable and every cyclic $S$- acts is locally cyclic.

\bpr  Every locally cyclic right $S$-act  is  hollow, and consequently, every cyclic right $S$-act is  hollow.\epr
\bpf Let $A_S$ be a locally cyclic $S$-act. If $A_S$ is simple, the result follows. Otherwise, Let $B$ be a proper subact of $A_S$. If  $C\cup B= A$ for some proper subact $C$ of $A$, take $a\in A\setminus B$ and $a'\in A\setminus C$. So there exists $a''\in A$ with $a,a'\in a''S$. Since $A= B\cup C$, we have $a''\in B$ or $a''\in C$ which implies that $a\in B$ or $a'\in C$, a contradiction. Thus  $C=A$, and $B$ is a superfluous subact of $A_S$.  \epf


\bte \lb{te000} A right $S$-act $A_S$ is  hollow if and only if $A_S$ is an indecomposable co-uniform right $S$-act.  \ete
\bpf Necessity. Suppose that $A_S$ is  hollow, and $B, C$ are proper subacts of $A$ such that $A=B\coprod C$. Thus $A=B\cup C$ which means that $B$ is not superfluous subact of $A$, a contradiction.

In view of Lemma \ref{le2.1}, the following the sufficiency is deduced.    \epf

In general being indecomposable does not imply being  hollow. For instance, Let $A_S$ be a cyclic $S$-act with a proper subact $B$, then $A\coprod^{B}A$ is indecomposable but not hollow. In particular, For a proper right ideal $I$ of a monoid $S$, $S\coprod^{I}S$ is indecomposable but not hollow. So we have the following strict implications.\\\centerline{cyclic  $\Longrightarrow$ locally cyclic $\Longrightarrow$  hollow  $\Longrightarrow$ indecomposable}

In the following proposition we characterize co-uniform $S$-acts.

\bpr  \lb{pr3.11} Every co-uniform $S$-act $A$ is indecomposable or $A=A_1\coprod A_2$, where each $A_i$ is simple.\epr
\bpf Suppose that $A_S$ is a co-uniform decomposable $S$-act. Let $A=\coprod_{i\in I} A_i$. If $|I|>2$, fix  $k\neq j\in I$ and put $B=A_k\coprod A_j$. So $B\cup ( \coprod_{i\neq j} A_i)=A$ and $B\cap ( \coprod_{i\neq j} A_i)=A_k\neq\emptyset$. Then $B$ is not coessential which is a contradiction. Thus $|I|=2$. Now, suppose that $A=A_1\coprod A_2$ such that $A_1$ is not simple. Let $B_1$ be a proper subact of $A_1$. Then $B= B_1\coprod A_2$ is a proper subact of $A$ such that $B\cap A_1\neq \emptyset$ and $B\cup A_1=A$ which means that $B$ is not coessential, a contradiction. Then $A=A_1\coprod A_2$ which $A_1, A_2$ are simple,  as desired.\epf

  Let $S$ be an arbitrary monoid and $A= \Theta\coprod\Theta\coprod\Theta$. Using Proposition \ref{pr3.11}, $A$ is not  co-uniform. So 
 for each arbitrary monoid $S$ there exists a finitely generated $S$-act which is not hollow or  co-uniform.

An $S$-act $A$ is said to be a \textit{uniserial} $S$-act if every two subacts of $A$ are comparable with respect to inclusion. In the next theorem we characterize an $S$-act  all its subacts are hollow.

\bte For an $S$-act $A_S$ the following statements are equivalent.\begin{itemize}
	\item[(i)] $A$ is a uniserial $S$-act.
	\item[(ii)] Every subact of $A$ is  hollow. 
	\item[(iii)] Every subact of $A$ generated by two elements is  hollow.\end{itemize}\ete
\bpf 

 The implications (i) $\lo$ (ii) and (ii) $\lo$ (iii) are obvious.

(iii) $\lo$ (i) Let $B$ and $C$ be subacts of $A$ and let $B\nsubseteq C$. Then there exists an element             $x\in B\backslash C$. To show that $C\subseteq B$,    suppose that $y\in C$.  Put  $N=xS\cup yS$. If $N=yS$, then $xS\subseteq N=yS\subseteq C$. So $x\in C$, a contradiction. Hence $yS$ is a proper subact of $N$,  and since $N$ is  hollow,   then $N=xS$. Therefore, $yS\subseteq N=xS\subset B$ which implies that $y\in B$, and so  $C\subseteq B$.\epf

\bpr \lb{pr3.111}  The following hold for a monoid $S$.
\begin{itemize}
	
	\item[(i)] Every indecomposable
	co-uniform $S$-act with a minimal generating set is cyclic.
	
	\item[(ii)] Every    hollow $S$-act with a minimal generating set is cyclic.
	\item[(iii)] Every     finitely generated hollow  $S$-act is cyclic.\end{itemize} \epr

\bpf It suffices to prove part (ii).   Let $A_{S}$ be a right $S$-act with a
minimal generating set $\{a_{i}\mid i\in I\}$. In contrary suppose that $|I|> 1$, and fix $i\in I$. Then  $a_iS\cup (\cup_{j\neq i} a_jS)= A$, and since $A_S$ is hollow, $A_S=\cup_{j\neq i} a_jS$, a contradiction.     \epf

Recall that a monoid $S$ satisfies condition (A) if all right
$S$-acts satisfy the ascending chain condition for cyclic subacts.
In \cite{kh} it is shown that a monoid  $S$ satisfies condition (A) if and only if  every locally cyclic $S$-act is cyclic, equivalently,  every right $S$-act contains a minimal generating set. Now, using this fact and the previous proposition we deduce the following result  as a generalization of that result in \cite{kh}.
\ble \lb{le4} A monoid $S$ satisfies condition (A) if and only if every  hollow $S$-act is cyclic. \ele 

We conclude this section considering the cover of  hollow $S$-acts. In \cite{kh}, it is shown that a cover of a locally cyclic right $S$-act is
indecomposable. Now, we extend this to the following result.
\ble \lb{le3} Each cover of a  hollow $S$-act is
indecomposable.\ele

\bpf Let 
$A_S$ be a  hollow $S$-act and $f:D_S\rightarrow A_S$  a coessential epimorphism. Suppose that $D= \coprod _{i\in
	I}D_i$ such that each $D_i$ is indecomposable. In contrary, suppose that  $|I| >1$ and choose
$i\neq j\in I$. Since $f\mid _{D\setminus D_i}$ is not an  epimorphism,  $ f(D\setminus D_i)$ is a proper subact of $A$ and 
$f(D\setminus D_i)\cup f(D\setminus D_j)=A$. Now since $A_S$ is  hollow,  $f(D\setminus D_j)=A$, and so $f\mid _{D\setminus D_j}$ is an epimorphism,
a contradiction. Therefore $B$ is indecomposable.\epf

The following corollary is a straightforward result of the
previous lemma.
\bco \lb{co1} For a monoid $S$ the following hold.

$~\,~(i)$ Every projective cover of a  hollow $S$-act is
cyclic.

$~\,(ii)$ Every strongly flat  (condition (P)) cover of a  hollow $S$-act is locally cyclic.\eco






{\tiny{\tiny{\tiny{\section{\large{\hspace{-0.7 cm}.\quad
			 The relation between	hollow	 and	Radical of  $S$-acts }}}}}}
In this section we consider local $S$-acts and the radical of an $S$-act. We also discuss the relationship between local and hollow $S$-acts. 
\bde A right $S$-act is called \textit{local} if it contains exactly one maximal subact. A monoid $S$ is also called \textit{right (left) local} if it contains exactly one maximal right (left) ideal.\ede

The set of maximal subacts of a right $S$-act $A_S$ is denoted by $\rm{Max}(A)$.
\ble \lb{le4.02} Every cyclic right $S$-act is simple or local.\ele
\bpf  Suppose that $A=aS$ is cyclic, and $A_S$ is not simple. By using Zorn's Lemma, $\rm{Max}(A) \neq \emptyset$. Now, suppose that $M\neq N$ are  maximal subacts of $A$.  Then  $M\cup N=A$ implies that  $a\in M$ or $a\in N$, and so $N=A$ or $M=A$, a contradiction. Thus $A$ is local.\epf

Now, we deduce the following remark which was discussed also in \cite{nak}.
\brm  Every monoid $S$ is  group or right local. Indeed the set $\{s\in S|~s~ \rm{is~ not~ right~ invertible}\}$ is either empty or the unique maximal right ideal of $S$.  Then right local and left local are equivalent for a monoid $S$. Thus we briefly call it local monoid.\erm



The following theorem  establishes a relation to   hollow $S$-acts with local and cyclic $S$-acts .
\bte Let $A_S$ be a right $S$-act. Then the following are equivalent:
\begin{itemize}

\item[(i)] $A_S$ is a  hollow right $S$-act and $\rm{Max}(A) \neq \emptyset$;
\item[(ii)]$A_S$ is a cyclic and local right $S$-act;
\item[(iii)] $A_S$ is a finitely generated local right $S$-act;
\item[(iv)] every  proper subact of $A_S$ is contained in a maximal subact, and $A_S$ is a local right $S$-act;
\item[(v)] $A_S$ contains a maximal subact $N$ such that $N\leq_s A$;
\item[(vi)] $A_S$ contains the maximum subact $N$ such that $N\leq_s A$.
\end{itemize}\ete
\bpf $(i) \lo (ii)$ Let $N$ be a maximal subact of $A_S$ and let $L$ be an arbitrary subact of
$A_S$ where $L \subsetneq N$. Since $N \cup L = A$, and $A_S$ is a  hollow right $S$-act, then $A = L$. Hence $A_S$ has
just one maximal subact. If $a\in A \setminus N$ and $L = aS$, then $A = aS$.

The implications $(ii) \lo (iii)$ and $(iii) \lo (iv)$  are obvious.

$(iv) \lo (v)$ Let $N$ be the unique maximal subact of $A$ and let $L$ be a proper subact of $A$. By assumption, $L\subseteq N$. Then $L\cup N=N\neq A$ and so $N\leq_s A$.

$(v)\lo (vi)$ Let $N$ be a maximal subact of $A$ which $N\leq_s A$ and let $B$ be a proper subact of $A$. So
$N \cup B \neq A$ and by maximality of $N$ we have $B\subseteq N$. So $N$ is maximum.

$(vi)\lo (i)$ Let $N$ be the maximum subact of $A$ which $N\leq_s A$. For each proper subact $B$ of $A$ we have $B\subseteq N \leq_s A$, we deduce that $B\leq_s A$. Therefore $A_S$ is  hollow.\epf

In general, every hollow (indecomposable co-uniform) $S$-act is not cyclic or local.
For instance, take  $S=(\mathbb{N}, \text{min})\cup
\{\varepsilon\}$ where $\varepsilon$ denotes the externally
adjoined identity greater than each natural element. Then $A=\{1,
2,3, ...\}$ is not cyclic act and $Max (A)=\emptyset$. But all its subacts are $\{1\}\subseteq\{1,2\}\subseteq\{1,2,3\}\subseteq...$, and so $A$ is  hollow.

Let $S$ be a monoid and $A$ a right $S$-act.  The radical of the act $A$ is the intersection of all maximal subacts of $A$,\\\centerline{   ${\displaystyle \mathrm {Rad} (A)=\cap \{N\mid N{\mbox{ is a maximal subact of A}}\}\,}$}. If $A$  contains no maximal subact, we put $\mathrm {Rad} (A)= A$. If $Rad(A)\neq \emptyset$, the $Rad(A)$ is a subact of $A$.

In module theory, radical submodule is equal to the union of superfluous submodules. The next proposition demonstrates that it is also valid for $S$-acts. To reach that we need the following lemma.

\ble If $a\in A$ and $C\leq A$ such that $aS\cup C=A$,      then $C=A $ or there exists a maximal subact $M$ of $A$ such that $C\subseteq M$ and $a\notin M$. \ele

\bpf Let $C\neq A$. Take $B=\{D|~D\lvertneqq A ~ and ~       C\subseteq D\}$. Clearly  $C\in B\neq\emptyset$  and    $B$   is a partially ordered set. 
Let $\{D_i\}_{i\in I}$ be a chain in $B$, so $D_i\lvertneqq A$
and $C\subseteq D_i$. Let $D=\cup_{i \in I} D_i$. If $D\lvertneqq A$, then $D$ is an upper bound. Otherwise, if $D=A$, $a\in A$ implies  $a\in D$, and there exists $ i\in I$ such that $a\in D_i$. Then $aS\subseteq D_i$ which implies that  $aS\cup D_i=D_i=A$, a  contradiction. Then   by  Zorn's Lemma, $B$ has a maximal element $M$.   So $M$ is a maximal subact of $A$ such that $C\subseteq M$,  $a\notin M$. \epf
 
 As we know, $A\leq_s A$ if and only if $A$ is simple.
\bpr \lb{pr3.33} Let $A_S$ be a right $S$-act. Then\\\centerline{
${\displaystyle \mathrm {Rad} (A)=\cup \{B\mid B\leq_s A\}\,}$.}\epr
\bpf Suppose that $\Gamma=\cup \{B\mid B\leq_s A\}$. First we show that $\Gamma\subseteq \mathrm {Rad} (A)$. If $\rm{Max}(A) = \emptyset$, clearly $\Gamma\subseteq \mathrm {Rad} (A)=A$. Otherwise, let $B\leq_s A$.   and $N$ be an arbitrary maximal subact of $A$. If $B\nsubseteq N$, being maximal of $N$ implies that $B\cup N=A$. Since $B\leq_s A$, $N=A$, a contradiction. Thus $B\subseteq N$, and so $\Gamma\subseteq \mathrm {Rad} (A)$. To show the converse, let $a\in Rad(M)$. First we show that $aS\leq_s A$.  If $aS=A$, then $A=Rad(A)$ and by Lemma \ref{le4.02} $A$ is simple. So  $aS=A\leq_s A$. Now, let $aS$ be a proper subact of $A$ and $aS\cup C=A$. If $C\neq A$ by previous lemma there exists a maximal subact $M$ of $A$ such that $C\subseteq M$ and $a\notin M$, but $a\in Rad(M)$  implies  $a\in M$,  a  contradiction. Then $C=A$ which means that $aS\leq_sA$. We deduce $ aS\subseteq \cup\{B|B\leq_sA\}$, and therefore $Rad(A)\subseteq\varGamma$.  \epf
Using the previous proposition, the following result is immediately deduced.

\bco For a monoid $S$  the following statements hold.
\begin{itemize}

\item[(i)] Let $A_S$ be a right $S$-act. Then for each element $a\in \mathrm{Rad}(A)$, $aS\leq_s A$. 
\item[(ii)]
  Let $A$ and $B$ be  right $S$-acts and let $f:A\lo B$ be an
$S$-monomorphism. Then $f(Rad(A)) \subseteq Rad(B)$.
\item[(iii)] $Rad(A)=A$ if and only if all finitely generated subact of $A$ are superfluous in $A$. 
\end{itemize} 
  \eco

  \bco   Let $A_S$ be a right $S$-act. Then each non-cyclic  hollow subact $B$ of $A$
  	is contained in $\mathrm {Rad} (A)$.\eco
  
  \bpf Assume that   $B$ is a  hollow subact of $A$ and $b\in B$.  So $bS$ is a proper subact of $B$ and $bS \leq_s B$, and by Lemma \ref{le2.2}, $bS\leq_s A$. Using the previous proposition, $bS\subseteq \mathrm {Rad} (A)$ which implies that $B\subseteq \mathrm {Rad} (A)$.  \epf
  
  Now, we give an equivalent condition for an $S$-act which its radical is superfluous.
\bte For a right $S$-act $A$ the following statements are equivalent. \begin{itemize}

\item[(i)]  $\mathrm {Rad} (A)\leq_s A$.

\item[(ii)] Every  proper subact of $A$ is contained in a maximal subact. \end{itemize}\ete

\bpf (i) $\lo$ (ii). Let  $C$ be a proper subact of $A$. Since $\mathrm {Rad} (A)\leq_s A$, $ \mathrm {Rad} (A)\cup C\neq A$. Suppose $\{M_i|~i\in I\}$ is the family of all maximal subacts of $A$. So $(\cap_{i\in I}M_i)\cup C\neq A$, which implies that $\cap_{i\in I}(M_i\cup C)\neq A$. Then there exists $j\in I$ such that $M_j\cup C\neq A$. Now, maximality of $M_j$ implies that $C\subseteq M_j$, and the result follows.

 (ii) $\lo$ (i). Suppose that $C$ is an arbitrary proper subact of $A$. There exists a maximal subact $M$ of $A$ with $C\subseteq M$. Then we have $C\cup \mathrm {Rad} (A)\subseteq M\cup \mathrm {Rad} (A)=M\neq A$. Thus, $\mathrm {Rad} (A)\leq_s A$.\epf

\bpr An $S$-act $A$ is finitely generated if and only if $A/Rad(A)$
 is finitely generated and $Rad(A)\leq_s A$.\epr

\bpf Let $A$ be finitely generated, clearly $A/Rad(A)$ is finitely generated. Let $C\leq A$, $Rad(A)\cup C=A$, by
Proposition \ref{pr3.33}, $Rad(A)=\cup\{B~|~B\leq_sA\}$, So $\cup\{B~|~B\leq_sA\}\cup C=A$. Since $A$ is finitely generated, there exist $B_1,...,B_m\leq_sA$ such that $B_1\cup B_2\cup ...B_m\cup C=A$. Since $B_1\leq_sA$ and $B_1\cup(B_2\cup...\cup B_m\cup C)=A$,  we imply that $B_2\cup...\cup B_m\cup C=A$. Since $B_2,...., B_m\leq_sA$, we continue this manner to imply $C=A$. Thus $Rad(A)\leq_sA$.

Sufficiency. Suppose that $A/Rad(A)=\cup_{i=1}^{i=n}[a_i]S$. So $Rad(A)\cup(\cup_{i=1}^{i=n}a_iS)=A$. Now, since $Rad(A)\leq_sA$, $\cup_{i=1}^{i=n}a_iS=A$. Thus $A$ is finitely generated.  \epf

{\tiny{\tiny{\tiny{\section{\large{\hspace{-0.7 cm}.\quad
						Supplemented Acts}}}}}}
In this section we introduce the notions of a supplement of a subact and supplemented $S$-acts, and general properties of them are discussed. Our  aim is to use the notion of a supplement of a subact to investigate the union of hollow $S$-acts.				
\bde Let $ B,C$  be proper subacts of a right $S$-act $A$. We call $C$ is a \textit{supplement} of $B$ in $A$, or $B$ has a supplement $C$ in $A$ if the following two conditions are satisfied.

\begin{itemize}
	
	\item[(i)]  $B\cup C=A$.
	\item[(ii)] If $D\subseteq C$ and $B\cup D=A$, then $D=C$.\end{itemize}
If every proper subact of $A$ has a supplement in $A$, then $A$ is called a \textit{supplemented} $S$-act.\ede

Clearly, If an $S$-act $A=B\coprod C$, then $C$ is a supplement of $B$. We first begin with elementary properties for  being supplement.
\ble\lb{le4.1} Let $A=B\cup C$. If $B\cap C\neq \emptyset$, Then $C$ is a supplement
of $B$ in $A$ if and only if  $C\cap B=\emptyset$ or $C\cap B\leq_s C$.\ele

\bpf Let $E$ be a subact of $C$. Then $(C\cap B)\cup E=C$ is equivalent to $A=B\cup E$ and so the result is easily checked. \epf
The following result presents that co-uniform implies supplemented.
\bpr \lb{th3.11} Every  co-uniform $S$-act is supplemented.\epr
\bpf Let $A$ be a right $S$-act and $B$ be a proper subact of
$A$. First suppose that $A$ is indecomposable. By Theorem \ref{te000}, $A$ is  hollow. Then $B\cup A=A$ and $(B\cap A)=B\leq_s A$ imply that
$A$ is a supplemented $S$-act. In the case that $A$ is not indecomposable, by Proposition \ref{pr3.11}, $A=B\coprod C$ where $B,C$ are simple acts. Thus $C$ is a supplement
of $B$.\epf

The converse of Proposition \ref{th3.11} is not valid. 
For instance,  let $S$ be an arbitrary monoid and $A= \Theta\coprod\Theta\coprod\Theta$. Using Proposition \ref{pr3.11}, $A$ is not  co-uniform. But, as all subsets of $A$ are also subacts, for each subact $B$ of $A$ we have $A\setminus B$  is a supplement of $B$. 

Let $C$ be a proper subact of an $S$-act $A$. By Lemma \ref{le2.2}, each superfluous subact of $C$ is also superfluous in $A$. So clearly $Rad(C)\subseteq C\cap Rad(A)$.
\bpr Suppose that  $C$ is a proper subact of an $S$-act $A$ such that $C$ is a supplement of a proper subact $B$ of $A$. Then the following hold.\begin{itemize}
	\item [(i)] If $D\cup C=A$ for some $D\subset B$, then $C$ is a supplement of $D$.
	\item [(ii)] If $A$ is finitely generated, then $C$ is also finitely generated.
	\item [(iii)] If $E$ is a subact of $C$ such that $E\leq_sA$, then $E\leq_sC$.
	\item [(iv)]  If $N\leq_sA$, then $N\cap C\leq_sC$.
	\item [(v)] If $N\leq_sA$, then $C$ is a supplement of $N\cup B$. 
	\item [(vi)] $Rad(C)=C\cap Rad(A)$.
\end{itemize}\epr 
\bpf (i) It is easily proved by using 
Lemmas \ref{le4.1} and \ref{le2.2}.

(ii)  Let $A$ be finitely generated. Since $B\cup C=A$, there is a finitely generated subact $X\subseteq C$ such that $B\cup X=A$. By the minimality of $C$, we imply that $C=X$.

(iii) Let $X$ be a subact of $ C$ with $E\cup X=C$.  Since $B\cup C=A$, we have $B\cup E\cup X=A$. Now, since $E\leq_sA$, $B\cup X=A$ and so  $X=C$.

(iv) Using part (iii) and Lemma \ref{le2.2}, it is clearly checked.

(v)  Let $N\leq_sA$. We have $(N\cup B)\cup C=A$. Let $X\subseteq C$ with $(N\cup B)\cup X=A$, then $N\leq_sA$ implies that $B\cup X=A$, and hence $X=A$.

(vi) We have $Rad(C)\subseteq C\cap Rad(A)$. To show the converse, if $N\leq_sA$,  by part (iv), $E=N\cap C\leq_sC$, and $E\subseteq Rad(C)$.  Therefore,\\  $C\cap Rad(A)= C\cap (\cup \{N\mid N\leq_s A\})=\cup \{N\cap C\mid N\leq_s A\} \subseteq Rad(C)$. 

\epf

Now, we turn our attention to the concept of supplement in a projective $S$-act.
\bpr Let $P$ be a projective $S$-act, and  $C$ be a supplement of $B$ in $P$. Then $C$ is projective or there exists an epimorphism $f:P\lo C$ such that $f(B)\leq_s C$.\epr

\bpf  Let $C$ be a supplement of $B$ in $P$. So $P=B\cup C$. If $B\cap C=\emptyset$, then $P=B\coprod C$, and $C$ is projective. Now, suppose that $B\cap C\neq\emptyset$. Let $\pi_1:C\lo C/(B\cap C)$ be the canonical epimorphism, and define $\pi_2:P\lo C/(B\cap C)$ by $ \pi_2(p)= \left\{
\begin{array}{ll} [p], & p\in C\\ \theta, & p\in B\end{array} \right.$. So since $P$ is projective, there exists a homomorphism $f:P\lo C$ with $\pi_1f=\pi_2$. It is easily checked that $\mathrm{Im}f\cup B=P$, and by assumption, $\mathrm{Im}f=C$. Moreover, since $f(B)\subseteq B\cap C\leq_s C$, by Lemma \ref{le2.2}, $f(B)\leq_s C$. \epf

Finally, we conclude this paper by considering the union of  hollow acts.
\bte Let $A$ be a right $S$-act such that $\mathrm {Rad} (A)\leq_s A$. The following statements are equivalent. \begin{itemize}

	\item[(i)] $A$ is a union of  hollow acts.

	\item[(ii)] Each proper subact $B$ of $A$ whose $A/B$ is finitely generated  has a supplement.
	\item[(iii)] Every maximal subact of $A$ has a supplement.  \end{itemize}\ete

\bpf (i) $\lo$ (ii). Suppose $A=\cup_{i\in I} L_i$ such that each $L_i$ is  hollow $S$-act. Let $B$ be a proper subact of $A$ such that $A/B$ is finitely generated. Then $A/B=\cup_{i\in I} (L_i\cup B)/B$. Since $A/B$ is finitely generated, $A=B\cup L_1\cup L_2\cup ... \cup L_n$ for some  hollow $S$-acts $L_1,L_2,...,L_n$ with $B\cap L_i\neq L_i$ for each $1\leq j\leq n$. Take $L=L_1\cup L_2\cup ... \cup L_n$. To show that $L$ is a supplement of $B$, let $X$ be a proper subact $L$. There exists $1\leq j\leq n$ such that $X\cap L_j$ is a proper subact of $L_j$. Now, since $L_j$ is  hollow, $(B\cap L_j)\cup (X\cap L_j)\neq L_j$. Thus $B\cup X\neq A$, and the result follows.

(ii) $\lo$ (iii)  follows by Lemma \ref{le1.1}. (iii) $\lo$ (i). Let $B$ be the union of all  hollow subacts of $A$. In contrary, suppose that $B$ is a proper subact of $A$. So there exists a maximal subact $N$ of $A$ with $B\subseteq N$. Let $L$ be a supplement of $N$ in $A$. If $L$ is simple, then $L\subseteq B$.  Otherwise, Let $X$ be a proper subact  of $L$. So $N\cup X\neq A$, and maximality of $N$ implies that $X$ is contained in $N$. So by Lemma \ref{le4.1},  $N\cap L\leq_s L$, and using Lemma \ref{le2.2}, $X\subseteq N\cap L \subseteq L$   implies  $X\leq_s L$. Then $L$ is a  hollow act. Therefore $L$ is contained in $B$, and so $A=L\cup N\subseteq B\cup N=N$, a contradiction. Therefore, $B=A$.
 Now suppose that $C$ is an arbitrary proper subact of $A$. There exists a maximal subact $M$ of $A$ with $C\subseteq M$. Then we have $C\cup \mathrm {Rad} (A)\subseteq M\cup \mathrm {Rad} (A)=M\neq A$. Thus, $\mathrm {Rad} (A)\leq_s A$.
\epf

\def\bibname{REFERENCES}


\begin{thebibliography}{99}
\bibitem{nak} Ahmadi, K., Madanshekaf, A.:  Nakayama's Lemma for acts over monoids. Semigroup Forum, {\bf{91}} (2), 321–337 (2015)
\bibitem{flu74}  Fluery, P.: Hollow modules and local endomorphism rings. Pacific J. Math. {\bf{53}}, 379-385 (1974)
 \bibitem{hara}  Harada, M.: On maximal submodules of a finite direct sum of  hollow modules III. Osaka J. Math. {\bf{22}},
81-95 (1985)
\bibitem{isb} Isbell, J.:  Perfect
monoids. Semigroup Forum. {\bf{2}}, 95-118 (1971)
\bibitem{kh} Khosravi, R., Ershad, M. and Sedaghatjoo, M.: Strongly flat and condition (P) covers of acts over
monoids. Comm. in Algebra.  {\bf{38}}, 4520-4530 (2010)  
\bibitem{kilp2000} Kilp, M. ,  Knauer, U. and Mikhalev, A.:  Monoids, Acts
and Categories. W. de gruyter. Berlin. (2000)

\bibitem{kilp96a} Kilp, M.: Perfect monoids revisited. Semigroup Forum. {\bf{53}},
225-229 (1996)

\bibitem{mah} Mahmoudi, M. and  Renshaw, J.: On covers of cyclic acts over
monoids. Semigroup Forum. {\bf{77}}, 325-338 (2008)
 
 \bibitem{uniform} Roueentan, M. and Sedaghatjoo, M.: On uniform acts over semigroups. 
  Semigroup Forum. {\bf{97}}, 229-243 (2018)
   \bibitem{wis}  Wisbauer, R.: Foundations of module and ring theory. Gordon and Breach Science Publishers, Reading (1991)
\end{thebibliography}
\end{document}